\documentclass{article}
\usepackage{amsfonts}
\usepackage{amsmath}
\usepackage{graphicx}
\usepackage{cite}
\usepackage{hyperref}
\usepackage{booktabs}

\begin{document}

\title{Asymmetric Duffing oscillator: jump manifold and border set}
\author{Jan Kyzio\l , Andrzej Okni\'{n}ski \\
Politechnika \'Swi\c{e}tokrzyska, Al. 1000--lecia PP 7, \\
25-314 Kielce, Poland}
\maketitle

\begin{abstract}
We study the jump phenomenon present in the forced asymmetric Duffing oscillator
using the known steady-state asymptotic solution. The major result is the 
computation of the jump manifold, which encodes global information about all
possible jumps.
\end{abstract}

\section{Introduction}
\label{intro}

In this work we study steady-state dynamics of the forced asymmetric Duffing
oscillator governed by the equation:%
\begin{equation}
\ddot{y}+2\zeta \dot{y}+\gamma y^{3}=F_{0}+F\cos \left( \Omega t\right) ,
\label{AsymDuffing}
\end{equation}%
where $\zeta $, $\gamma $, $F_{0}$, $F$ are paremeters and $\Omega $ is the
angular frequency of the periodic force, which has a single equilibrium
position and a corresponding one-well potential \cite{Kovacic2011}. This
dynamical system in particular and Duffing-type equations in general, which
can be used to describe pendulums, vibration absorbers, beams, cables,
micromechanical structures, and electrical circuits, have a long history 
\cite{KovacicBrennan2011}. The equation of motion (\ref{AsymDuffing}) can
describe several non-linear phenomena, such as various non-linear resonances,
symmetry breaking, chaotic dynamics, period-doubling route to chaos,
multistability and fractal dependence on initial conditions, and jumps \cite%
{Hayashi1964,Ueda1980,Ueda1979,Szemplinska1986,KovacicBrennan2011}.

The aim of our work is to research the jump phenomenon described and 
investigated for the system (\ref{AsymDuffing}) in an interesting study by
Kovacic and Brennan \cite{Kovacic2011}, using the implicit function
machinery. An early approach to the jump phenomenon is due to Holmes and Rand 
\cite{Holmes1976} in the Catastrophe Theory setting where a discriminant of
a cubic equation, arising in the amplitude-frequency response function, was
analysed to detect multivaluedness of the amplitude. Recently, Kalm\'{a}r-Nagy 
and Balachandran applied a differential condition to detect vertical
tangencies, characteristic for the jump phenomenon \cite{Kalmar2011}.

In the next Section we describe the steady state solution (\ref{steady-state}%
) \cite{Szemplinska1986,Jordan1999,Kovacic2011} which is an implicit
function of $A_{0}$ and $\Omega $ (\ref{A0}) and compute an implicit
function of $A_{1}$ and $\Omega $ (\ref{A1-a}), (\ref{A1-b}). 
Working in the framework developed in our earlier papers \cite{Kyziol2021,Kyziol2022} we compute the
jump manifold in Section \ref{jumps} containing information about all
possible jumps which is the main achievement of this work. We summarize our
results in the last Section.

\section{The steady-state solution}
\label{steady}

The steady-state solution of Eq. (\ref{AsymDuffing}) of form:%
\begin{equation}
y\left( t\right) =A_{0}+A_{1}\cos \left( \Omega t+\theta \right) ,
\label{steady-state}
\end{equation}%
can be computed by any of asymptotic methods to yield non-linear algebraic
equations expressing variables $A_{0},A_{1},\theta $ and $\Omega $ by
parameters $\zeta $, $\gamma $, $F_{0}$, $F$ \cite%
{Szemplinska1986,Jordan1999,Kovacic2011}: 
\begin{subequations}
\label{1a1b1c}
\begin{eqnarray}
-A_{1}\Omega ^{2}+3\gamma A_{0}^{2}A_{1}+\frac{3}{4}\gamma A_{1}^{3}-F\cos
\theta &=&0,  \label{1a} \\
-2\zeta A_{1}\Omega -F\sin \theta &=&0,  \label{1b} \\
\gamma A_{0}^{3}+\frac{3}{2}\gamma A_{0}A_{1}^{2}-F_{0} &=&0.  \label{1c}
\end{eqnarray}
\end{subequations}
Eliminating $\theta $ from Eqs. (\ref{1a}), (\ref{1b}) we get two implicit
equations for $A_{0},\ A_{1}$ and $\Omega $:
\begin{subequations}
\label{1d1e}
\begin{eqnarray}
A_{1}^{2}\left( -\Omega ^{2}+3\gamma A_{0}^{2}+\frac{3}{4}\gamma
A_{1}^{2}\right) ^{2}+4\Omega ^{2}\zeta ^{2}A_{1}^{2} &=&F^{2},  \label{1d}
\\
\gamma A_{0}^{3}+\frac{3}{2}\gamma A_{0}A_{1}^{2}-F_{0} &=&0.  \label{1e}
\end{eqnarray}
\end{subequations}
Computing $A_{1}^{2}$ from Eq. (\ref{1e}) and substituting into (\ref{1d})
we obtain an implicit equation for $A_{0},\ \Omega $:
\begin{equation}
f\left( \Omega ,A_{0};\gamma ,\zeta ,F,F_{0}\right)
=\sum\nolimits_{k=0}^{9}c_{k}A_{0}^{k}=0,  \label{A0}
\end{equation}%
where coefficients $c_{k}$ are given in Table \ref{tab:T1}.
\begin{table}[ht]
\caption{Coefficients of polynomial (\ref{A0}}
\label{tab:T1}
\centering
\begin{tabular}{|l|l|}
\hline
$c_{9}=25\gamma ^{3}$ & $c_{4}=16\Omega ^{2}\gamma F_{0}$ \\ \hline
$c_{8}=0$ & $c_{3}=-9\gamma F_{0}^{2}+6\gamma F^{2}$ \\ \hline
$c_{7}=-20\Omega ^{2}\gamma ^{2}$ & $c_{2}=-4F_{0}\Omega ^{4}-16\zeta
^{2}\Omega ^{2}F_{0}$ \\ \hline
$c_{6}=-15\gamma ^{2}F_{0}$ & $c_{1}=4\Omega ^{2}F_{0}^{2}$ \\ \hline
$c_{6}=-15\gamma ^{2}F_{0}$ & $c_{0}=-F_{0}^{3}$ \\ \hline
\end{tabular}
\end{table}

We can also obtain implicit equation for $A_{1},\ \Omega $. We solve the
cubic equation Eq. (\ref{1e}) for $A_{0}$ computing one real root (two other
roots are complex): 
\begin{subequations}
\begin{equation}
A_{0}=-\frac{A_{1}^{2}}{2Y}+Y,\quad Y\equiv \sqrt[3]{\sqrt{\frac{1}{8}%
A_{1}^{6}+\frac{1}{4\gamma ^{2}}F_{0}^{2}}+\frac{1}{2\gamma }F_{0}}.
\label{A1-a}
\end{equation}%
Then we substitute $A_{0}$ from Eq. (\ref{A1-a}) into Eq. (\ref{1d}),
obtainig finally a complicated but useful implicit equation for $A_{1},\
\Omega $: 
\begin{equation}
g\left( \Omega ,A_{1};\gamma ,\zeta ,F,F_{0}\right) =A_{1}^{2}\left( 3\gamma
A_{0}^{2}+\tfrac{3}{4}\gamma A_{1}^{2}-\Omega ^{2}\right) ^{2}+4\Omega
^{2}\zeta ^{2}A_{1}^{2}-F^{2}=0,  \label{A1-b}
\end{equation}%
where $A_{0}$ and $Y$ are defined in (\ref{A1-a}).

\section{Singular points}
\label{singular}

Conditions for singular points are \cite{Kyziol2021,Kyziol2022}:

\end{subequations}
\begin{subequations}
\label{3a3b3c}
\begin{eqnarray}
f\left( \Omega ,A_{0};\gamma ,\zeta ,F,F_{0}\right) &=&0,  \label{3a} \\
\frac{\partial f\left( \Omega ,A_{0};\gamma ,\zeta ,F,F_{0}\right) }{%
\partial \Omega } &=&0,  \label{3b} \\
\frac{\partial f\left( \Omega ,A_{0};\gamma ,\zeta ,F,F_{0}\right) }{%
\partial A_{0}} &=&0.  \label{3c}
\end{eqnarray}%
Since the derivative $\frac{\partial }{\partial \Omega }f\left( \Omega
,A_{0},\gamma ,\zeta ,F,F_{0}\right) $ factorizes:

\end{subequations}
\begin{equation}
\frac{\partial f\left( \Omega ,A_{0};\gamma ,\zeta ,F,F_{0}\right) }{%
\partial \Omega }=8A_{0}\Omega \left( F_{0}-\gamma A_{0}^{3}\right) \left(
F_{0}+A_{0}\left( 5\gamma A_{0}^{2}-4\delta ^{2}-2\Omega ^{2}\right) \right)
\label{factor}
\end{equation}
we get from Eq. (\ref{3b}):

\begin{equation}
F_{0}=2\Omega ^{2}A_{0}+4\zeta ^{2}A_{0}-5\gamma A_{0}^{3},  \label{F_0}
\end{equation}%
other possibilities leading to trivial solutions or cases without a solution.

Substituting (\ref{F_0}) into Eqs. (\ref{3a}), (\ref{3c}) we obtain: 
\begin{subequations}
\label{XA_0}
\begin{gather}
\begin{array}{l}
256\zeta ^{4}X^{4}+1280\zeta ^{6}X^{3}+\left( -96c\zeta ^{2}+2304\zeta
^{8}\right) X^{2} \\ 
+\left( -336c\zeta ^{4}+1792\zeta ^{10}\right) X+512\zeta ^{12}-240\zeta
^{6}c+45c^{2}=0,%
\end{array}
\label{X} \\
\frac{1}{45F^{2}\gamma ^{2}}A_{0}^{2}=16\zeta ^{2}X^{3}+64\zeta
^{4}X^{2}+\left( 9\gamma F^{2}+80\zeta ^{6}\right) X+15\gamma F^{2}\zeta
^{2}+32\zeta ^{8}.  \label{A_0}
\end{gather}%
where $X=\Omega ^{2}$, $c=\gamma F^{2}$.

We have checked, plotting the implicit function of $X$, $\zeta $, $c$
defined by Eq. (\ref{X}), that Eq. (\ref{X}) has no solutions $X=\Omega
^{2}>0$. We thus conclude that the system of equations (\ref{1d1e}) has no
singular points.

\section{The jump phenomenon}
\label{jumps}

\subsection{Jump conditions and jump manifold}
\label{manifold}

Jump conditions in implicit function setting read  \cite{Kyziol2021,Kyziol2022}:

\end{subequations}
\begin{subequations}
\label{J1J2}
\begin{eqnarray}
f\left( \Omega ,A_{0};\gamma ,\zeta ,F,F_{0}\right) &=&0,  \label{J1} \\
\frac{\partial f\left( \Omega ,A_{0};\gamma ,\zeta ,F,F_{0}\right) }{%
\partial A_{0}} &=&0.  \label{J2}
\end{eqnarray}
\end{subequations}
where equation (\ref{J2}) is the condition for a vertical tangency (note
that equations (\ref{J1J2}) are Eqs. (\ref{3a}), (\ref{3c})).

Solving equations (\ref{J1J2}) we obtain:

\begin{subequations}
\label{SOL}
\begin{gather}
J\left( A_{0};\gamma ,\zeta ,F,F_{0}\right)
=\sum\nolimits_{k=0}^{21}a_{k}A_{0}^{k}=0,  \label{Sol1} \\
\Omega ^{2}=\tfrac{\left( -50\gamma ^{4}\right) A_{0}^{12}+95\gamma
^{3}F_{0}A_{0}^{9}+\left( 6F^{2}\gamma ^{2}-39\gamma ^{2}F_{0}^{2}\right)
A_{0}^{6}+\left( 3F^{2}\gamma F_{0}-7\gamma F_{0}^{3}\right)
A_{0}^{3}+F_{0}^{4}}{2A_{0}\left( F_{0}-10\gamma A_{0}^{3}\right) \left(
F_{0}-\gamma A_{0}^{3}\right) ^{2}},  \label{Sol2}
\end{gather}
\end{subequations}

the non-zero coefficients $a_{k}$ of the polynomial $J\left(
A_{0}\right) $\ given in Table \ref{tab:T2}.

\begin{table}[ht]
\centering
\caption{Non-zero coefficients of polynomial (\ref{Sol1})}
\label{tab:T2}
\begin{tabular}{|l|l|}
\hline
$a_{21}=4000\gamma ^{7}\zeta ^{2}$ & $a_{9}=3248\gamma ^{3}\zeta
^{2}F_{0}^{4}-72F^{2}\gamma ^{3}\zeta ^{2}F_{0}^{2}$ \\ \hline
$a_{18}=-16\,000\gamma ^{6}\zeta ^{2}F_{0}$ & $a_{8}=36F^{4}\gamma
^{3}F_{0}-978F^{2}\gamma ^{3}F_{0}^{3}$ \\ \hline
$a_{17}=600F^{2}\gamma ^{6}$ & $a_{6}=528\gamma ^{2}\zeta
^{2}F_{0}^{5}-240F^{2}\gamma ^{2}\zeta ^{2}F_{0}^{3}$ \\ \hline
$a_{15}=23\,880\gamma ^{5}\zeta ^{2}F_{0}^{2}-480F^{2}\gamma ^{5}\zeta ^{2}$
& $a_{5}=9F^{4}\gamma ^{2}F_{0}^{2}+138F^{2}\gamma ^{2}F_{0}^{4}$ \\ \hline
$a_{14}=-1920F^{2}\gamma ^{5}F_{0}$ & $a_{3}=24F^{2}\gamma \zeta
^{2}F_{0}^{4}-152\gamma \zeta ^{2}F_{0}^{6}$ \\ \hline
$a_{12}=768F^{2}\gamma ^{4}\zeta ^{2}F_{0}-15\,512\gamma ^{4}\zeta
^{2}F_{0}^{3}$ & $a_{2}=-6F^{2}\gamma F_{0}^{5}$ \\ \hline
$a_{11}=36F^{4}\gamma ^{4}+2166F^{2}\gamma ^{4}F_{0}^{2}$ & $a_{0}=8\zeta
^{2}F_{0}^{7}$ \\ \hline
\end{tabular}
\end{table}

The polynomial $J\left( A_{0}\right) $, complicated as it is, encodes global
information about all possible jumps. We shall thus refer to equation (\ref%
{Sol1}), which defines an implicit function of variables $A_{0}$, $\gamma $, 
$\zeta $, $F$, $F_{0}$, as jump manifold equation. Thus, the jump manifold $%
\mathcal{J}\left( A_{0},\gamma ,\zeta ,F,F_{0}\right) $:
\begin{equation}
\mathcal{J}\left( A_{0},\gamma ,\zeta ,F,F_{0}\right) \ \mathcal{=\ }\left\{
\left( A_{0},\gamma ,\zeta ,F,F_{0}\right) :\quad J\left( A_{0};\gamma
,\zeta ,F,F_{0}\right) =0\right\} ,  \label{Jm}
\end{equation}%
belongs to $5D$ space. It is purposeful to introduce projection of the jump
manifold onto the parameter space:%
\begin{equation}
\mathcal{J}_{\perp }\left( \gamma ,\zeta ,F,F_{0}\right) \ \mathcal{=}\text{ 
}\left\{ \left( \gamma ,\zeta ,F,F_{0}\right) : \text{there is such 
}A_{0}\text{ that }J\left( A_{0};\gamma ,\zeta ,F,F_{0}\right) =0\right\} .
\label{Jm_projected}
\end{equation}%
In other words, for any set of parameters $\gamma ,\zeta ,F,F_{0}$ belonging
to $\mathcal{J}_{\perp }$ there is a jump in dynamical system (\ref%
{AsymDuffing}) and all jumps occur for parameters belonging to $\mathcal{J}%
_{\perp }$.

We shall consider $2D$ and $3D$ projections, plotting $\mathcal{J}\left(
A_{0};\gamma _{\ast },\zeta _{\ast },F_{\ast },F_{0}\right) $ and $\mathcal{J%
}\left( A_{0};\gamma _{\ast },\zeta _{\ast },F,F_{0}\right) $, respectively,
where parameters $\gamma _{\ast },\zeta _{\ast },F_{\ast }$ or $\gamma
_{\ast },\zeta _{\ast }$\ are fixed.

\subsubsection{$2D$ projection, $J\left( A_{0};\protect\gamma _{\ast },%
\protect\zeta _{\ast },F_{\ast },F_{0}\right) =0$}
\label{2D}

\begin{figure}[h!]
\center
\includegraphics[width=10.5cm, height=7cm]{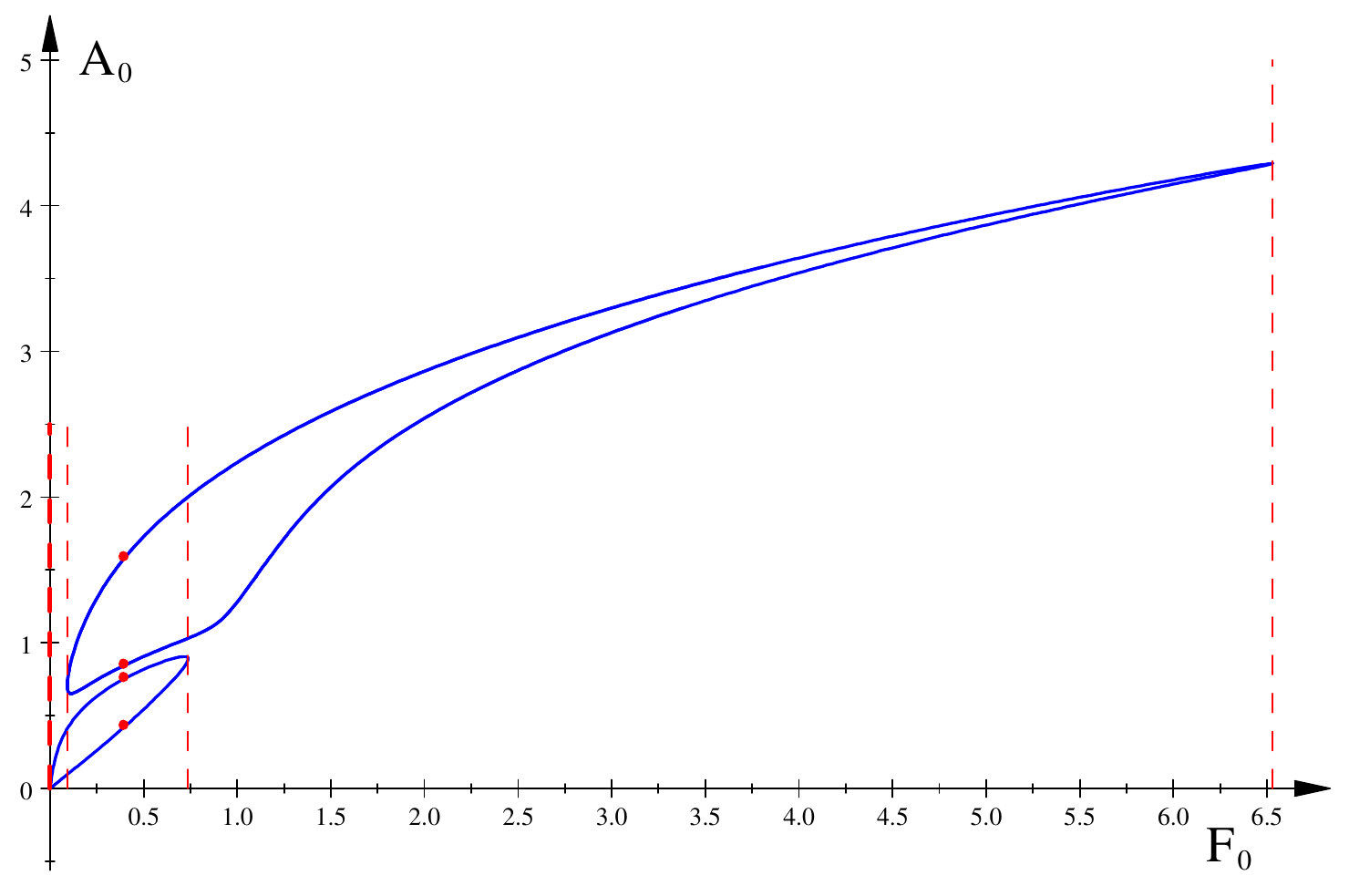} 
\caption{Jump manifold $\mathcal{J}\left( A_{0};\gamma _{\ast },\zeta _{\ast
},F_{\ast },F_{0}\right) $, $\gamma _{\ast }=0.0783$, $\zeta _{\ast }=0.025$%
, $F_{\ast }=0.1$ (blue) and four border points on intersections of $%
\mathcal{J}$ and vertical red lines.
}
\label{F1}
\end{figure}

Global picture of the jump manifold $\mathcal{J}\left( A_{0};\gamma
_{\ast },\zeta _{\ast },F_{\ast },F_{0}\right) $ where $\gamma _{\ast
}=0.0783$, $\zeta _{\ast }=0.025$, $F_{\ast }=0.1$ and $A_{0}$, $F_{0}$ are
variable is shown in Fig. \ref{F1}.

Firstly, all points lying on the blue curve correspond to jumps. Moreover,
there are four critical points, dividing Fig. $1$ into parts and referred to
as \textit{border points}: $F_{0}^{\left( 1\right) }=0,\ F_{0}^{\left(
2\right) }=0.0920,\ F_{0}^{\left( 3\right) }=0.7385,\ F_{0}^{\left( 4\right)
}=6.\,5321$, defined and computed in the Subsection \ref{border}. It turns
out that $F_{0}\in \left( F_{0}^{\left( 1\right) },\ F_{0}^{\left( 2\right)
}\right) $ there are two vertical tangencies, for $F_{0}\in \left(
F_{0}^{\left( 2\right) },\ F_{0}^{\left( 3\right) }\right) $ there are four,
for $F_{0}\in \left( F_{0}^{\left( 3\right) },\ F_{0}^{\left( 4\right)
}\right) $ there are two and there are no vertical tangencies for $%
F_{0}>F_{0}^{\left( 4\right) }$.

\begin{figure}[h!]
\center
\includegraphics[width=9cm, height=6cm]{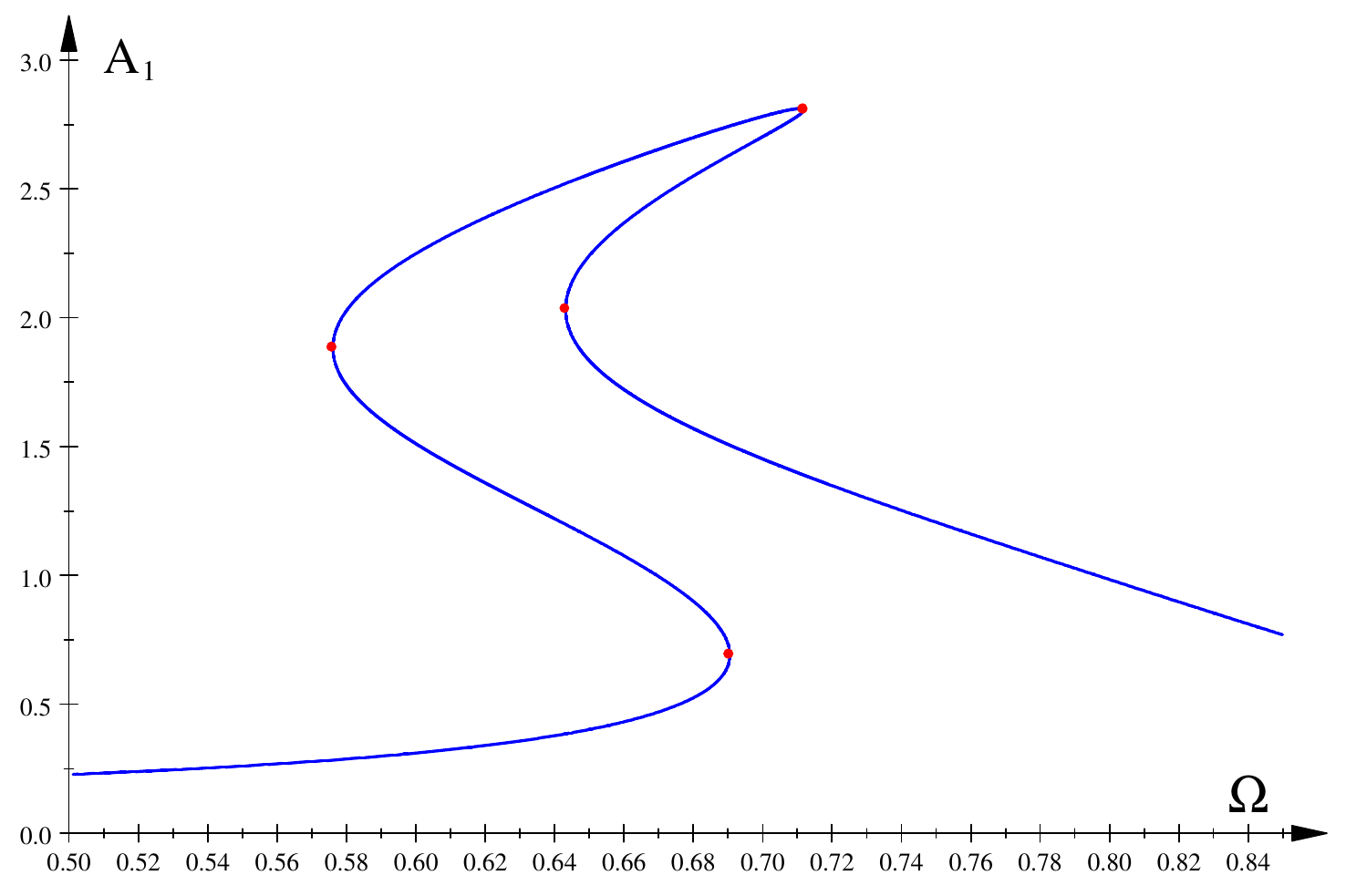}
\caption{Amplitude-frequency response curve $A_{1}\left( \Omega \right) $, $\gamma =0.0783$, $\zeta =0.025$, $F=0.1$%
, $F_{0}=0.4$.
}
\label{F3}
\end{figure}

For example, in Fig. \ref{F3} the case $F_{0}=0.4$ is shown. 
More exactly,
the implicit function $A_{1}\left( \Omega \right) $, computed with help of
Eq. (\ref{A1-b}) is plotted for $\gamma =0.0783$, $\zeta =0.025$, $F=0.1$
and $F_{0}=0.4$. Red dots, denoting vertical tangencies, correspond to red
dots in Fig. \ref{F1}. These points can be easily computed from Eqs. (\ref%
{J1J2}), (\ref{1d1e}). 

Indeed, solving equations (\ref{J1J2}) for $\gamma =0.0783$, $\zeta =0.025$, 
$F=0.1$, $F_{0}=0.4$ we get four real solutions $\Omega $, $A_{0}$ shown in the 
first two column in Table \ref{tab:T3}. Then for the above values of $\Omega 
$ we solve equations (\ref{1d1e}) obtaining four values of $A_{1}$ listed in
the third column of the Table \ref{tab:T3}. 

\newpage

\begin{table}[ht]
\centering
\caption{Solutions of Eqs. (\ref{J1J2}) and (\ref{1d1e})}
\label{tab:T3}
\begin{tabular}{|c|c|c|}
\hline
$\Omega $ & $A_{0}$ & $A_{1}$ \\ \hline\hline
\multicolumn{1}{|l|}{$0.576\,122\,891$} & \multicolumn{1}{|l|}{$%
0.846\,633\,527$} & \multicolumn{1}{|l|}{$1.\,882\,759\,746$} \\ \hline
\multicolumn{1}{|l|}{$0.643\,209\,846$} & \multicolumn{1}{|l|}{$%
0.755\,260\,872$} & \multicolumn{1}{|l|}{$2.\,032\,001\,367$} \\ \hline
\multicolumn{1}{|l|}{$0.690\,545\,624$} & \multicolumn{1}{|l|}{$%
1.\,583\,776\,750$} & \multicolumn{1}{|l|}{$0.691\,474\,188$} \\ \hline
\multicolumn{1}{|l|}{$0.711\,882\,658$} & \multicolumn{1}{|l|}{$%
0.425\,889\,574$} & \multicolumn{1}{|l|}{$2.\,806\,379\,023$} \\ \hline
\end{tabular}
\end{table}

\begin{figure}[h!]
\center
\includegraphics[width=9cm, height=6cm]{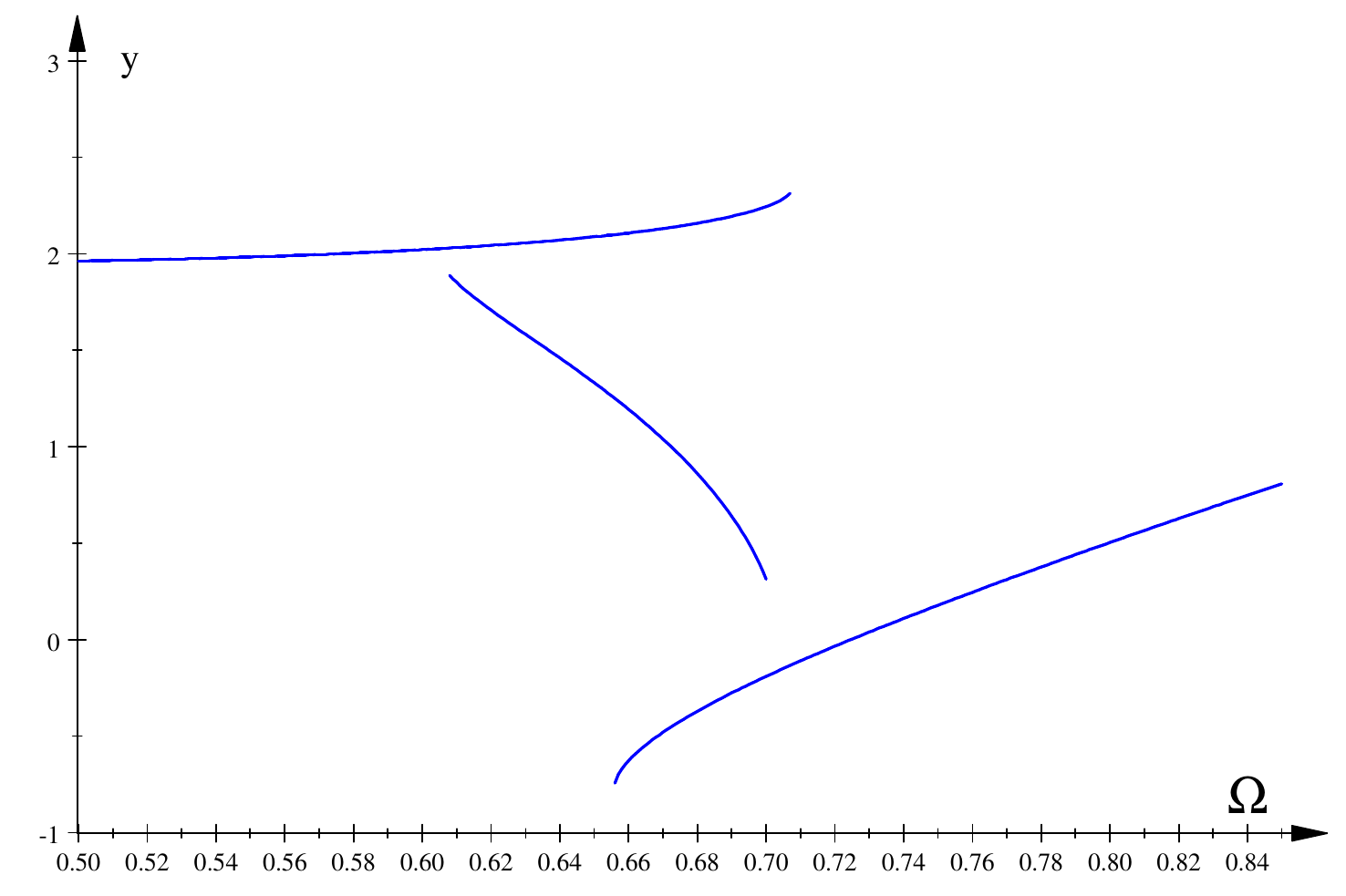}
\caption{Bifurcation diagram, $\gamma =0.0783$, $\zeta =0.025$, $F=0.1$%
, $F_{0}=0.4$. }
\label{F10}
\end{figure}

In Fig. \ref{F10} the bifurcation diagram is shown for the set of parameters listed in Fig. \ref{F3}. 
Note good agreement of end points of bifurcation branches in Fig. \ref{F10} with $\Omega$ 
coordinates of red dots in Fig. \ref{F3}.

\subsubsection{$3D$ projection, $J\left( A_{0};\protect\gamma _{\ast }, 
\protect\zeta _{\ast },F,F_{0}\right) =0$}
\label{3D}

\begin{figure}[h!]
\center
\includegraphics[width=12cm, height=8cm]{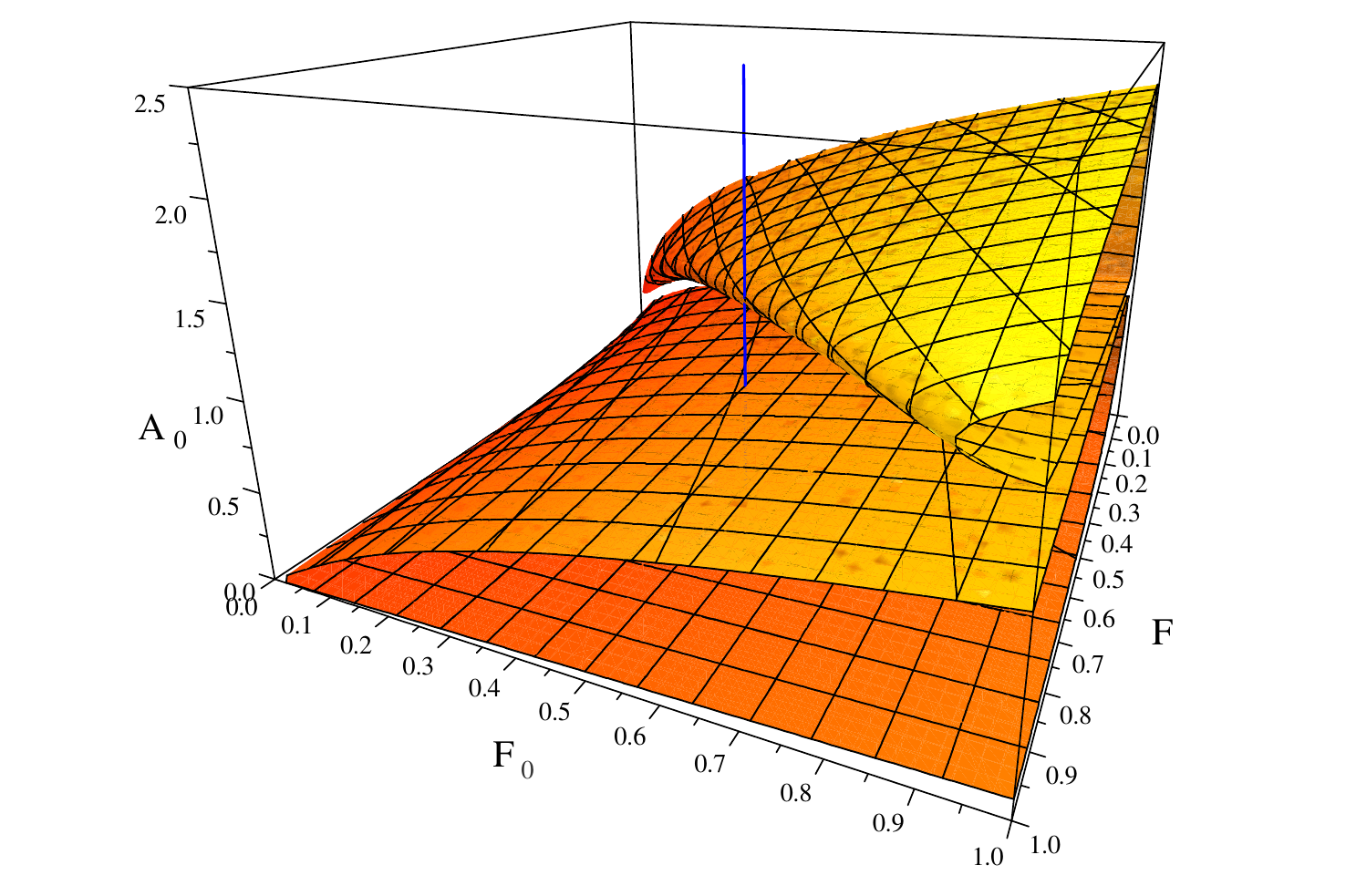} \vspace{-0.6cm}
\caption{Jump manifold $\mathcal{J}\left( A_{0};\gamma _{\ast },\zeta _{\ast
},F,F_{0}\right) $, $\gamma _{\ast }=0.0783$, $\zeta _{\ast }=0.025$.
}
\label{F4}
\end{figure}

We now fix two parameters only, for example $\gamma _{\ast }=0.0783$, $%
\delta _{\ast }=0.025$, and plot the jump manifold $\mathcal{J}\left(
A_{0};\gamma _{\ast },\zeta _{\ast },F,F_{0}\right) $ as a $3D$ surface, see Fig. \ref{F4}.

Next we compute one border point. For example, we choose $F_{0}=0.5$ $\left(
\gamma _{\ast }=0.0783,\delta _{\ast }=0.025\right) $ and compute the
corresponding border point as $F=0.544\,860$, $A_{0}=1.\,238\,340$ as
explained in the next Subsection.

The blue vertical line, given by the just computed data $\left( 0.544\,860,\ 0.5,\
A_{0}\right) $ with $A_{0}$ variable, touches the upper lobe of the jump
manifold exactly at the border point $\left( 0.544\,860,\ 0.5,\
1.\,238\,340\right) $.

\subsection{Border sets}
\label{border}

We shall now determine condition for the border set: the set of points in
the parameter space $\left( \gamma ,\ \zeta ,\ F,\ F_{0}\right) $, such that
number of vertical tangencies changes at these points. Condition for the
border set is that the polynomial $J\left( A_{0}\right) $ given in Eq. (\ref%
{Sol1}) and Table $2$ has multiple roots. Qualitative behaviour of the
polynomial equation $J\left( A_{0}\right) $ can be seen in Figs. \ref{F1}, \ref{F3} 
where $2D$ and $3D$ projections of the implicit function $J\left(
A_{0};\gamma ,\zeta ,F,F_{0}\right) =0$ are shown. To find parameters'
values for which the polynomial $J\left( A_{0};\gamma ,\zeta ,F,F_{0}\right) 
$ has multiple roots we demand that resultant of $J\left( A_{0}\right) $ and
its derivative $J^{\prime }\left( A_{0}\right) =\frac{d}{dA_{0}}J\left(
A_{0}\right) =\sum\nolimits_{k=0}^{20}b_{k}A_{0}^{k}$ is zero \cite{Gelfand2008,Janson2010}: 
\begin{equation}
R\left( J,J^{\prime };\gamma ,\zeta ,F,F_{0}\right) =0.  \label{R1}
\end{equation}

Resultant of the polynomial (\ref{Sol1}), Table $2$ is a determinant of the $%
\left( m+n\right) \times \left( m+n\right) $ Sylwester matrix, $n=21$, $m=20$%
:%
\begin{equation}
R\left( J,J^{\prime };\gamma ,\zeta ,F,F_{0}\right) =\det \left( 
\begin{array}{ccccccc}
a_{n} & a_{n-1} & a_{n-2} & \ldots & 0 & 0 & 0 \\ 
0 & a_{n} & a_{n-1} & \ldots & 0 & 0 & 0 \\ 
\vdots & \vdots & \vdots &  & \vdots & \vdots & \vdots \\ 
0 & 0 & 0 & \ldots & a_{1} & a_{0} & 0 \\ 
0 & 0 & 0 & \ldots & a_{2} & a_{1} & a_{0} \\ 
b_{m} & b_{m-1} & b_{m-2} & \ldots & 0 & 0 & 0 \\ 
0 & b_{m} & b_{m-1} & \ldots & 0 & 0 & 0 \\ 
\vdots & \vdots & \vdots &  & \vdots & \vdots & \vdots \\ 
0 & 0 & 0 & \ldots & b_{1} & b_{0} & 0 \\ 
0 & 0 & 0 & \ldots & b_{2} & b_{1} & b_{0}%
\end{array}%
\right)  \label{R2}
\end{equation}%
and is an enormously complicated polynomial in variables $\gamma ,\ \zeta ,\
F_{0},\ F$. However, if we fix three paramaters, say $\gamma ,\ \zeta ,\ F$,
then the equation $R\left( J,J^{\prime }\right) =0$ can be solved
numerically and thus critical values of $F_{0}$ can be computed.

For example, we have solved equation (\ref{R1}), $R\left( J,J^{\prime
};\gamma _{\ast },\zeta _{\ast },F_{\ast },F_{0}\right) =0$, for 

\begin{figure}[h!]
\center
\includegraphics[width=6cm, height=4cm]{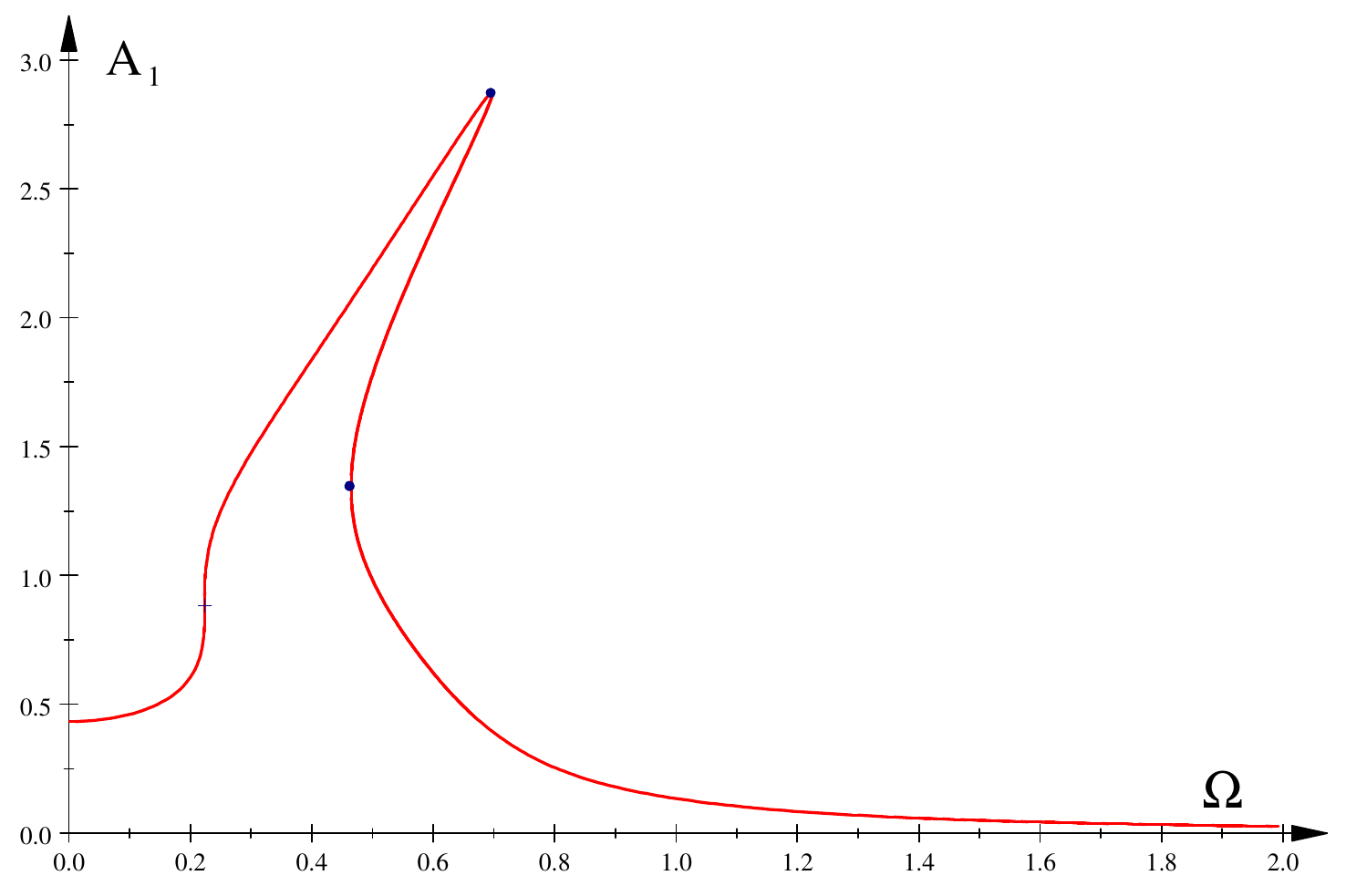}
\includegraphics[width=6cm, height=4cm]{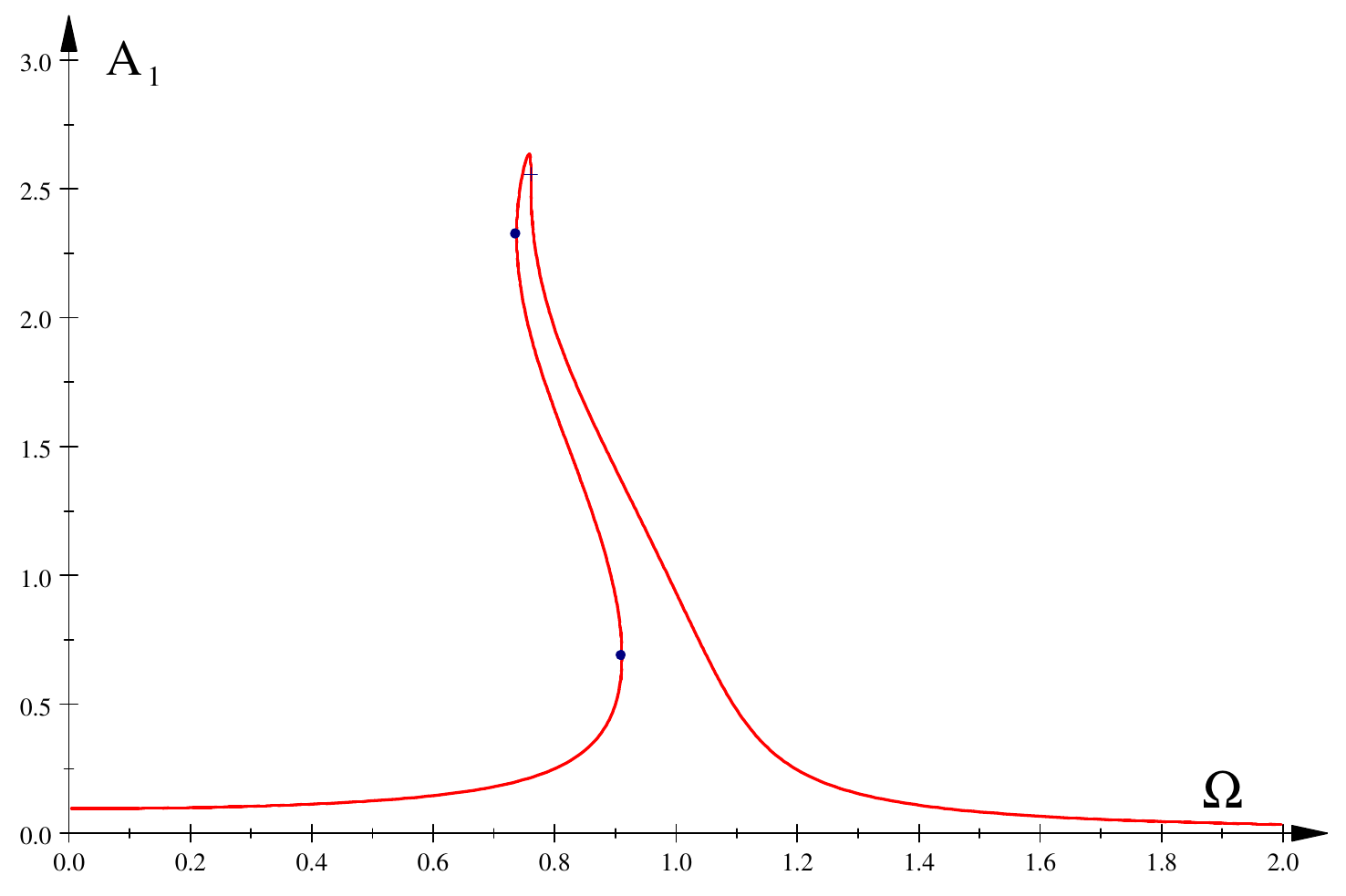}
\includegraphics[width=6cm, height=4cm]{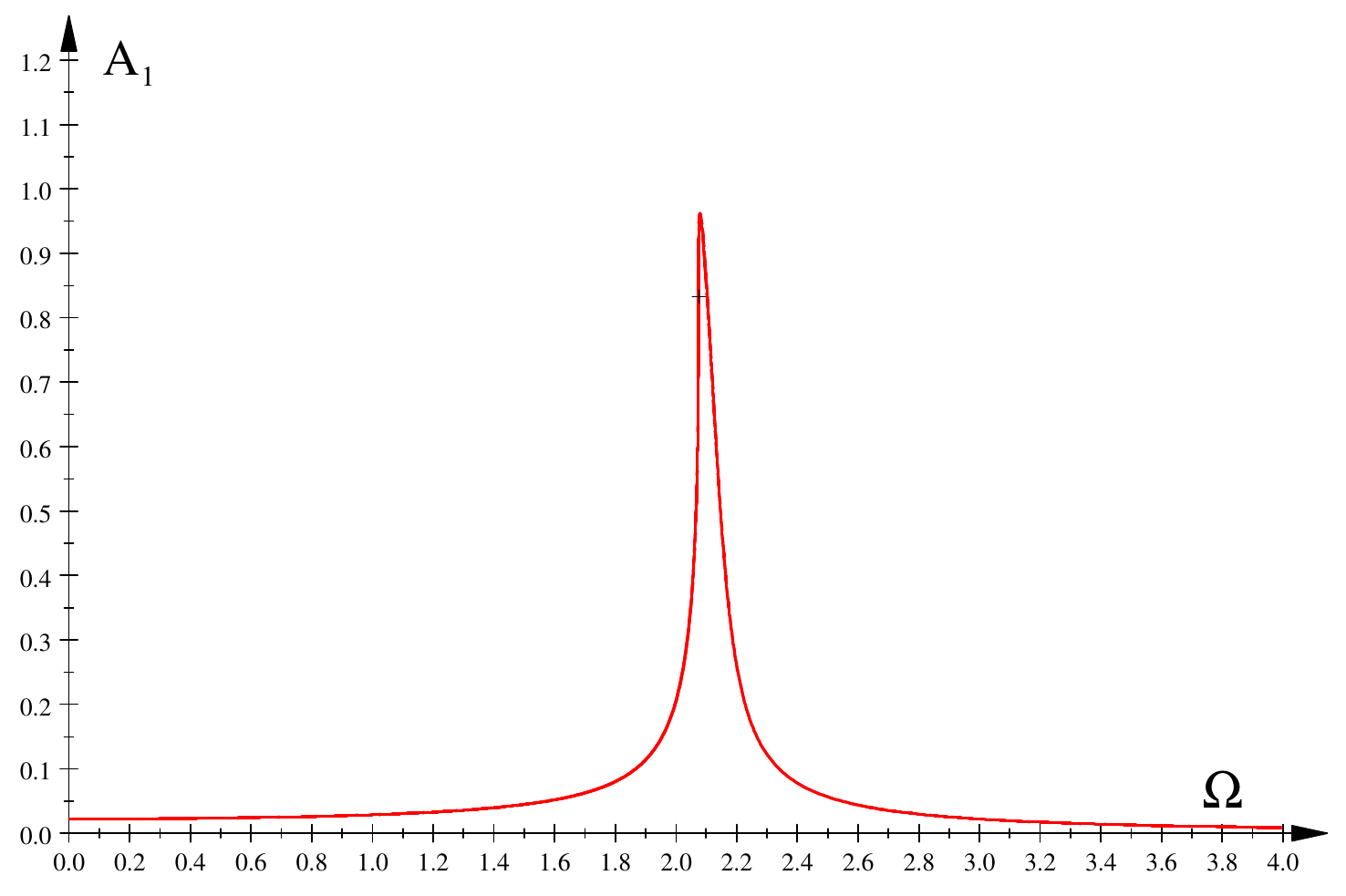}
\caption{Amplitude-frequency response curves: 
 $\gamma =0.0783$, $\zeta =0.025$, $F=0.1$, 
 $F_{0}^{\left( 2\right) }=0.092$ (left top), 
 $F_{0}^{\left( 3\right) }=0.7385$ (right top), 
 $F_{0}^{\left( 4\right) }=6.532$ (bottom).}
\label{F5}
\end{figure}

\noindent $\gamma_{\ast }=0.0783$, $\zeta _{\ast }=0.025$, $F_{\ast }=0.1$ obtaining the
following real positive solutions:
$F_{0}^{\left( 1\right) }=0$, $%
F_{0}^{\left( 2\right) }=0.092\,075$, $F_{0}^{\left( 3\right) }=0.738\,510$, 
$F_{0}^{\left( 4\right) }=6.\,532\,092$. 
In Figs. \ref{F5}  border amplitudes for 
$\gamma =0.0783$, $\zeta =0.025$, $F=0.1$ and $F_{0}^{\left( 2\right) }$, 
$F_{0}^{\left( 3\right) }$, $F_{0}^{\left( 4\right) }$ are marked with blue crosses, respectively.
Blue dots denote points of jumps.

We have also solved equation (\ref{R1}), $R\left( J,J^{\prime };\gamma
_{\ast },\zeta _{\ast },F,F_{0\ast }\right) =0$, for $\gamma _{\ast }=0.0783$%
, $\zeta _{\ast }=0.025$, $F_{0\ast }=0.5$ obtaining real positive
solutions: $F^{\left( 1\right) }=0$, $F^{\left( 2\right) }=0.026\,998\,9$, $%
F^{\left( 3\right) }=0.077\,925\,6$, $F^{\left( 4\right) }=0.544\,859\,5$.
Next, for $\gamma =0.0783$, $\zeta =0.025$, $F_{0}=0.5$, and $F=0.544\,859\,5
$ we have computed from Eqs. (\ref{J1J2}) the border value $A_{0}=1.\,238\,340$, see
the blue vertical line in Fig. \ref{F4}.

\subsection{Number of solutions of Eq. (\protect\ref{A0}) for a given value
of $\Omega $}
\label{number}

There are also other qualitative changes of the amplitudes $A_{1}\left(
\Omega \right) $ controlled by the parameters. For example, number of
solutions of Eq. (\ref{A0}) for a given value of $\Omega $ may change. This
happens when two vector tangencies appear at the same value of $\Omega $.
For example, let $\gamma =0.0783$, $\zeta =0.025$, $F=0.1$. To find a value
of $F_{0}$ for which this happens we have to find a double root of $\Omega $
of equations (\ref{J1J2}). Therefore, solving Eqs. (\ref{J1J2}) numerically
for several values of $F_{0}$\ we easily find that for $F_{0}=0.301\,007$
there is indeed a double root: $\Omega =0.597\,114$, $A_{0}=0.679\,284$ and $%
\Omega =0.597\,114$, $A_{0}=1.\,411\,787$. There is another such case, for $%
F_{0}=0.429\,166$ there is a double root: $\Omega =0.714\,419$, $%
A_{0}=0.459\,118$ and $\Omega =0.714\,419$, $A_{0}=1.\,628\,271$.

\begin{figure}[h!]
\center
\includegraphics[width=6cm, height=4cm]{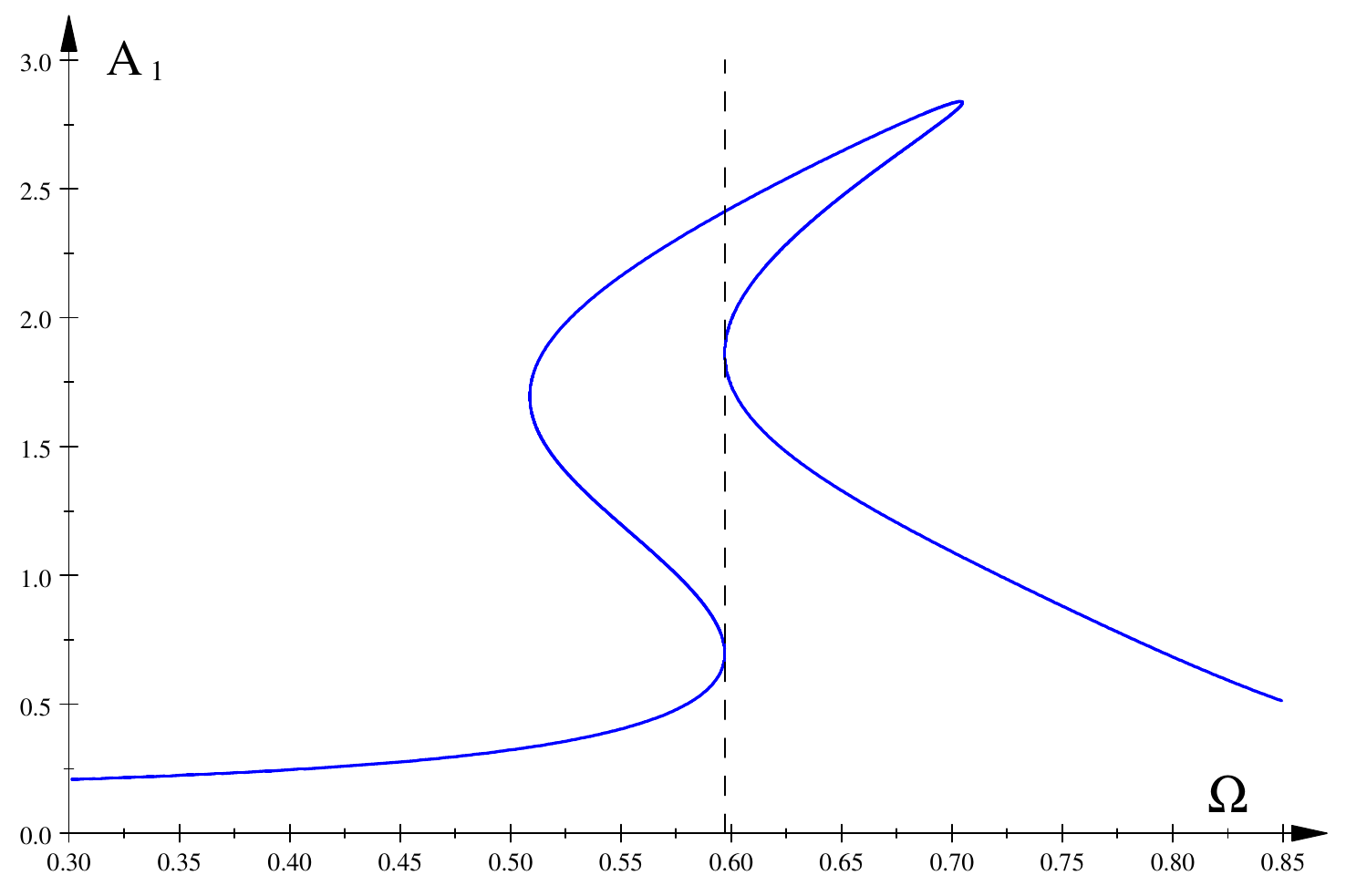}
\includegraphics[width=6cm, height=4cm]{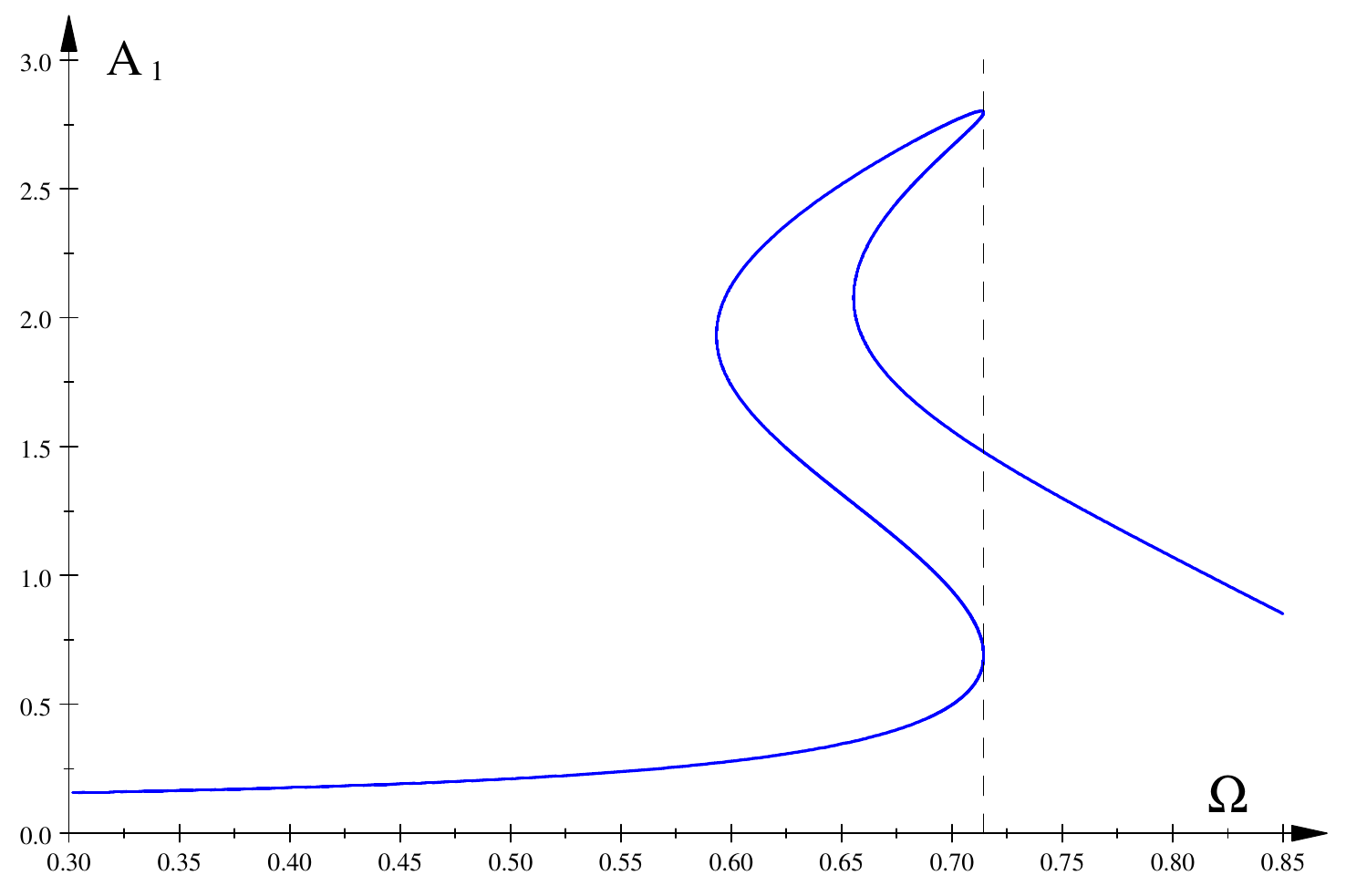}
\caption{Amplitude-frequency response curves $A_{1}\left( \Omega \right) $: $\gamma
=0.0783$, $\zeta =0.025$, $F=0.1$, $F_{0}=0.301$ (left), $F_{0}=0.429$
(right).}
\label{F8}
\end{figure}

\noindent Therefore, for $F_{0}\in \left( 0.301,\ 0.429\right) $ equation (\ref{A0}) has five solutions.


\section{Summary}
\label{summary}

Working in the implicit function framework \cite{Kyziol2021,Kyziol2022},  
we have computed in Subsection \ref{manifold}, using the (approximate) steady-state solution obtained in 
\cite
{Szemplinska1986,Jordan1999,Kovacic2011},  the jump manifold \ref{Jm},  comprising information about all jumps in 
the dynamical system \ref{AsymDuffing}.  Our formalism, described in Section \ref{jumps} 
 -- built on an idea to use a differential 
condition to detect vertical tangencies, due to Kalm\'{a}r-Nagy and Balachandran \cite{Kalmar2011} --
 can be applied to arbitrary steady-state solution. 

Our work on the asymmetric Duffing oscillator is  a supplementation 
and amplification of results obtained by Kovacic and Brennan \cite{Kovacic2011}. More precisely, 
the sequence of figures 8.4 (a) --(e), computed in  \cite{Kovacic2011} for $\gamma =0.0783$, $\zeta =0.025$, $F=0.1$ and 
$F_{0}=0.01$,$\ 0.2$,$\ 0.4$, $0.5$, $0.95$, respectively, can be appended with Figs. \ref{F5} and \ref{F8}, computed 
for $F_{0}=0.092$, $0.7385$, $6.532$ and $F_{0}=0.301$, $0.429$. Therefore, the whole sequence of 
metamorphoses of the curve $A_{1}\left( \Omega \right) $ consists of plots computed for 
$F_{0}=\mathbf{0.01}$,$\ 0.092$, $\mathbf{0.2}$,$\ 0.301$, $\mathbf{0.4}$, $%
0.429$, $\mathbf{0.5}$, $0.7385$, $\mathbf{0.95}$, $6.532$ where numbers highlighted in bold 
correspond to Figs. 8.4 (a) -- (e) plotted in \cite{Kovacic2011}.


\end{document}